\documentclass[12pt]{paper}

\usepackage{srcltx}
\usepackage{amsmath}
\usepackage{amssymb}
\usepackage{color}

\def\be#1\ee{\begin{equation}#1\end{equation}}
\newcommand{\bea}{\begin{eqnarray}}
\newcommand{\eea}{\end{eqnarray}}
\newcommand{\beas}{\begin{eqnarray*}}
\newcommand{\eeas}{\end{eqnarray*}}

\newtheorem{lemma}{Lemma}
\newtheorem{proposition}{Proposition}

\newtheorem{conjecture}{Conjecture}

\def\square{\hfill\hbox{\vrule\vbox{\hrule\phantom{o}\hrule}\vrule}}
\newenvironment{proof}{\textit{Proof.}}{\qquad\square\vspace{2mm}}

\begin{document}

\title{On generalized discrete  PML optimized  for propagative and evanescent  waves}

\author{Vladimir Druskin, \thanks{Schlumberger-Doll Research}
\and
Murthy Guddati  \thanks{North Carolina
State University}
\and
Thomas Hagstrom,
\thanks{Southern Methodist University}}
\maketitle
\begin{abstract}
We suggest a  unified  spectrally matched optimal grid approach   for finite-difference and finite-element approximation of the PML.
The new approach allows to combine optimal discrete absorption  for  both evanescent and propagative waves. 
\end{abstract}

\section{Introduction}
We approximate the Neumann-to-Dirichlet (NtD) map of  wave problem
in unbounded domain. After Fourier transform we obtain,
    \be\label{eq1} u_{xx}- \lambda u=0, \qquad x\in [0,\infty] \ee
and due to the infinity condition we are limited to outgoing wave
solutions,
    \[u=ce^{-\sqrt{\lambda}x}\]
satisfying NtD condition,
    \be\label{NtD}\frac{u}{u_x}|_{x=0}=-\frac{1}{\sqrt\lambda}.\ee
Here $\lambda=\kappa^2-\omega^2,$ where $\kappa$ and $\omega$ are
respectively (tangential) spacial and temporal frequencies. Also,
(\ref{eq1}) can be equivalently rewritten in the first order form
as,
    \be\label{1st} u_x=sv, \ v_x=su,\ee
where $s=\sqrt{\lambda}$. In terms of $u$ and $v$,  condition
(\ref{NtD}) can be equivalently rewritten as,
    \be\label{NtD1st}   \frac{u}{v}|_{x=0}=-1.\ee

The NtD can be numerically realized via rational approximation
theory using several approaches \cite[etc]{Lind75, EM79, IDK, GT00,
AsDrGuKn, GuLi06, ZaGu06, Hagstrom}. In \cite{IDK, AsDrGuKn} and
\cite{GuLi06, ZaGu06} this approximant was realized as respectively
finite-difference (FD) and  finite-element (FE) discretization of an
absorbing layer similar to well known Perfectly Matched Layer (PML)
\cite{Berenger}. In particular, the FD scheme was designed as an
optimal rational approximant separately for evanescent  solutions
corresponding to $\lambda\ge 0$ \cite{IDK} and propagative waves
\cite{AsDrGuKn} corresponding to $\lambda<0$, but not for the both
types of the solutions simultaneously. On the other hand, FE
approach is more flexible; while \cite{GuLi06} focuses on
propagative waves, it was shown in \cite{ZaGu06} that both
propagative and evanescent waves can be treated simultaneously. Most
recently, these FD and FE approximations are interpreted as special
quadrature rules with complete wavefield approximation...
\cite{Hagstrom}.

In this paper, we show that simultaneous treatment of propagative
and evanescent waves is possible not only in FE setting, but also in
FD setting. The key to this observation is a recently-discovered
equivalence between the FE and FD approaches for the two-sided
problem. Utilizing this link, we present two alternative approaches
to implement the NtD map and comment on their relative merits.
Furthermore, utilizing Zolotorev approximation theory and complete
wavefield approximation interpretation, we present an NtD map that
is an optimal approximation for propagating as well as evanescent
waves.

The outline of the paper is as follows. We start in section 2 with
the overview of optimal rational approximation of the NtD map by
considering both propagative and evanescent waves. Section 3
contains the description of FE and FD approximations of two-sided
problems and the equivalence between them. In section 4, we consider
rational approximation of the NtD map of the exterior problem and
present FE and FD realizations. The implementation details in time
domain and relative merits of the (or three) approaches are
considered in section 5. Numerical examples are presented in section
6. Finally, section 7 concludes the paper with some closing remarks.
{\bf (where do we fit complete wavefield approximation of Tom?)}

\section{Optimal Rational Approximation of NtD map for Propagative and Evanescent Waves}

Let us for simplicity consider time-harmonic case with $\omega=1$
and consider time-dependent problems later. Let us present our
rational approximant of $-\lambda^{-1/2}$ as
    \be\label{p/q} -\lambda^{-1/2} \approx
            R(\lambda)=p(\lambda)/q(\lambda), \ee
where $p$ and $q$ are polynomials of degrees $K-1$ and $K$
respectively. Introducing a polynomial of degree $N=2K$ given by,
    \[h(s)=sp(s^2)+q(s^2),\]
with $s=\sqrt{\lambda}$, we transform (\ref{p/q}) to Newman function
  \be\label{Newman}R(s^2)=p(s^2)/q(s^2)=\frac{h(s)-h(-s)}{s[h(s)+h(-s)]}
            = \frac{-1+h(s)/h(-s)}{s[1+h(s)/h(-s)]}.\ee
Then the  relative error of the NtD map is approximately
proportional to the reflection coefficient,
    \[\frac{h(s)}{h(-s)}.\]
According to (\ref{NtD1st}), the exact solution $(u,v)$ of
(\ref{1st}) is proportional to $(1,-1)$. In reality, due to the
approximation error, $(u,v)|_0=c_1(1,-1)+c_2(1,1)$, and
$\frac{c_2}{c_1}=\frac{h(s)}{h(-s)}$, i.e., the reflection
coefficient is the ratio of the incoming and outgoing waves.

Minimization of $\frac{h(s)}{h(-s)}$ on a real positive  interval is
the classical first Zolotarev problem solved in 1872. Zolotarev's
solution was first applied to the optimal FD approximation of the
NtD map for evanescent solutions in \cite{IDK} and then to the
approximation of propagative modes in \cite{AsDrGuKn}. The ABC for
both propagative and evanescent waves should approximate the true
NtD map on both negative $[-1,\lambda_1]$ and positive
$[\lambda_2,\lambda_3]$ intervals. They respectively correspond to
intervals $S_p=[\sqrt{-1}, \sqrt{\lambda_1}]$ and
$[S_e=\sqrt{\lambda_2},\sqrt{\lambda_3}]$ of variable $s$. The so
called spectrally matched finite-difference scheme (a.k.a FD
Gaussian spectral rule or optimal FD grid) \cite[etc]{DrKn} allows
arbitrary $h(s)$, but does not simultaneously treat propagative and
evanescent waves. On the other hand, propagative and evanescent
waves have been simultaneously treated using FE approximation in
\cite{ZaGu06}. The specific approximation is based on linear FE
approximation with
 midpoint integration \cite{GuLi06}, which is linked to special rational approximation \cite{GT00} with
    \be\label{droot}h(s)=t(s)^2,\ee
where $t$ is a polynomial of degree $k$.\footnote{However, as it
will be shown in the Section~4, simultaneous treatment of
propagating and evanescent waves is even possible with more general
$h(s)$, if the FD approach is used.}  Hence, considering the success
in \cite{ZaGu06}, we limit the current treatment to the restricted
form of $h$ in (\ref{droot}). With such restriction, minimization of
$\max_{s\in S_e\cup S_p}\left|\frac{h(s)}{h(-s)}\right|$ is
equivalent to solving
    \be\label{zolc}\min_{\deg t\le k} \max_{s\in S_e\cup S_p} \left|\frac {t(s)}{t(-s)}\right |.\ee

It is well known that the necessary and sufficient conditions for
optimality of a real rational approximant on a real interval is
so-called the Equal Ripple Theorem (ERT) \cite{PP}. It says that the
optimal error of [(K-1)/K] approximant has $2k-1$ zeros and $2k$
equal absolute value alternating extrema on the interval of
optimality. Generally, there is no similar result for complex
rational approximation \cite{Varga}. Here, instead of minimizing
(\ref{zolc}), we construct an approximant based on  classical
Zolotarev  results. We hope that its error is close to(\ref{zolc}).

If $t=t_et_p$, $\deg t_e=l<k$, $\deg t_p=k-l$, where $t_e$ and $t_p$
have respectively (non-coinciding) roots on $S_e$ and $S_p$, then
$\left| \frac {t(s)}{t(-s)}\right |$ has $2k+1$ maxima on $S_p\cup
S_e$. Moreover, $\left| \frac {t_e(s)}{t_e(-s)} \right| = 1$ on
$S_p$ and $\left| \frac  {t_p(s)}{t_p(-s)} \right| = 1$ on $S_e$,
which implies that,
    \[\max_{S_e}\left|\frac {t(s)}{t(-s)}\right |= \max_{S_e}\left|\frac {t_e(s)}{t_e(-s)}\right |, \]
and
    \[\max_{S_p}\left|\frac {t(s)}{t(-s)}\right |=\max_{S_p}\left|\frac {t_p(s)}{t_p(-s)}\right | .\]

Thus, we can take as $t_e$ and $t_p$ as the classical optimal
Zolotarev approximants on $S_e$ and $S_p$ respectively, and obtain
the quality of the total approximation the same as the one of the
separate problems.

The remaining question is: can the constructed approximant be
optimal in global sense, or, at least, how close is its error to
(\ref{zolc}). Obviously, $\max_{S_p}\left|\frac {t(s)}{t(-s)}\right
|$ and   $\max_{S_e}\left|\frac {t(s)}{t(-s)}\right |$ may be
different. Varying  $l$ one can equate $\max_{S_p}\left|\frac
{t(s)}{t(-s)}\right |$ and $\max_{S_e}\left|\frac
{t(s)}{t(-s)}\right |$  for a countable set of arrays
$\lambda_1,\lambda_2,\lambda_3$.

Here we conjecture that the ERT can be extended to the first
Zolotarev problem on two intervals in $C$ in the following way.
\begin{conjecture} Let $t_e(s)/t_e(-s)$ and $t_p(s)/t_p(-s)$ be the solutions of the Zolotarev problems on $S_e$ and $S_p$ respectively.

1. There are infinitely many arrays $\lambda_1,\lambda_2, \lambda_3$
for which there exists $l$, such that
    \be\label{czol}\max_{S_e}\left|\frac {t_e(s)}{t_e(-s)}\right |=\max_{S_p}\left|\frac {t_p(s)}{t_p(-s)}\right |.\ee

2. If (\ref{czol}) is valid, then  $t=t_et_p$ solves (\ref{zolc}).
\end{conjecture}
Results of \cite{leonid's reference} indicate that, if (\ref{czol})
is valid, then at
 least the approximant is optimal in the Cauchy--Hadamard  sense. Generally, it is always possible to find $l$ such that $\max_{S_p}\left|\frac {t(s)}{t(-s)}\right |$ and $\max_{S_e}\left|\frac {t(s)}{t(-s)}\right |$ are of the same order, in which case, it is natural to assume that the approximation error will be of the order of (\ref{zolc}).

\section{Equivalence of FE and FD Approximations for Two-sided Problems}

While the emphasis of this paper is on the approximation of the
one-sided problem on  $[0,\infty)$, in this section, we consider the
two-sided problem on $[0,1]$ and show that there exist equivalence
between spectrally matched FD grids and midpoint integrated linear
FE mesh. We then utilize these results in Section 4 to construct an
effective NtD map for the one-sided problem on $[0,\infty)$.

\subsection{Continuum problem}

{\bf QUESTION: You have used * for many row vectors and matrices.
Should we be just using transpose?} {\color{blue} decide later}.

Let us consider eq. (\ref{1st}) on $[0,1]$, and define the two-sided
DtN map as matrix-valued function $F(s)\in C^{2\times 2}$
    \[F(s)u_b=v_b,\]
where $u_b=[u(0),u(1)]^*$, $v_b=[v(0),-v(1)]^*$. It is easy to see
that $(u,v)$ is a linear combination of,
    \[(e^{\pm sx},\pm e^{\pm sx}),\]
and simple computation shows that,
    \be\label{impex} F(s)= \frac {1}{\sinh(s)}\left[\begin{array}{cc}\cosh(s) & -1 \\ -1&\cosh(s)\end{array}\right]=
                Z\left[\begin{array}{cc}\tanh(s/2) & 0 \\ 0 & \coth(s/2)\end{array}\right]Z^* ,\ee
where $Z$ is an orthogonal matrix
    \[Z=\frac{1}{\sqrt 2}\left[\begin{array}{cc}1 &- 1 \\ 1 & 1 \end{array}\right].\]
Similarly, we define  propagator operator from left to right as
matrix-valued function $G(s) \in C^{2\times 2}$ $G(s)w(0)=w(1)$,
where $w=(u,v)^*$ and from (\ref{impex}) we obtain
    \be\label{prop} G
        = \left[\begin{array}{cc} \cosh(s) & \sinh(s) \\ \sinh(s) & \cosh(s) \end{array} \right]
        = Z\left[\begin{array}{cc}\exp(s) & 0 \\ 0 & \exp(-s) \end{array} \right] Z^*.
    \ee

\subsection{Discrete problem: linear FE mesh with midpoint rule}

It was shown in \cite{GuLi06} that the discretization of the
original second-order from in (\ref{eq1}) with midpoint-integrated
linear FE mesh would lead to exponential convergence of the NtD map.
Furthermore, it was shown in \cite{GuDr07} that such a FE
discretization is equivalent to Crank-Nicholson discretization of
the first order form (\ref{1st}), i.e.
    \be\label{1stfe} \frac{u_{i+1}-u_i}{l_i}
        = s\frac{v_{i+1}+v_i}{2}, \ \frac{v_{i+1}-v_i}{l_i}
        = s\frac{u_{i+1}+u_i}{2}, \qquad i=1,\ldots n.
    \ee
where $l_i$, $i=1,\ldots, n$ are the FE lengths with
$\sum_{i=1}^nl_i = 1$. It can be easily verified that $(u_j,v_j)$,
$j=1,\ldots,n$, is a linear combination of
    \[\left(\prod_{i=1}^j\frac{1\pm l_is/2}{1\mp l_is/2},
            \pm \prod_{i=1}^j\frac{1\pm l_is/2}{1\mp l_is/2}\right) .\]
Comparing the above (approximate) solution with the exact solution
and noting that $\sum_{i=1}^nl_i = 1$, the FE solution approximates
the exponential as
 \[\exp(s)\approx exp(s)= t(-s)/t(s),\]
 where
    \[t(s)=\prod_{i=1}^n(1- l_is/2),\]
Assuming
    \be\label{feimp} u(0) = u_1, \ u(1)=u_n,\quad
        v(0)\approx v_1, \ v(1)\approx v_n,\ee
we can compute the approximate NtD map  as,
    \be\label{impfe} \tilde F(s)
        = \frac {1}{sinh(s)} \left[\begin{array}{cc}cosh(s) & -1 \\ -1 & cosh(s)\end{array}\right]
        = Z\left[\begin{array}{cc}tanh(s/2) & 0 \\ 0 & coth(s/2)\end{array}\right]Z^* .
    \ee
Here,
    \[sinh(s)=\frac{exp(s)- exp(-s)} {2}\approx \sinh(s),\qquad
        cosh(s)=\frac{exp(s)+ exp(-s)} {2}\approx \cosh(s), \]
    \[tanh(s/2)=\frac{t(s)-t(-s)}{t(s)+t(-s)}\approx \tanh(s/2), \qquad
        coth(s)=1/tanh(s)\approx \coth(s).\]
Similarly, the discrete propagator from left to right matrix can be
computed as
    \be\label{propfe} \tilde G
        = \left[\begin{array}{cc}cosh(s) & sinh(s) \\ sinh(s) & cosh(s)\end{array}\right]
        = Z\left[\begin{array}{cc}exp(s) & 0 \\ 0 & exp(-s)\end{array}\right]Z^* .\ee
Vectors $\frac{1}{\sqrt 2}(1,\pm 1)$ are the eigenvectors of $\tilde
G$, so it has so called fixed point property, i.e., if
$u(0)/v(0)=\pm 1$ then $u(1)/v(1)=\pm 1$ and vice versa. This
implies that, if exact half-space BC (\ref{NtD1st}) is applied at
$x=0$,  it will be also valid at $x=1$ regardless of the accuracy of
the FE approximation. In other words, adding an FE-discretized
interval to a half-space does not alter the NtD map of the
half-space. Furthermore, it was shown in \cite{GuLi06} that adding a
midpoint-integrated finite element to an approximate half-space can
only decrease the approximation error in the NtD map. This property
was used in \cite{ZaGu06} to enhance the approximation, originally
designed for propagative waves, to simultaneously absorb evanescent
waves.

\subsection{Discrete Problem: spectrally matched finite-difference grids}
It was shown in \cite{} that one-sided, two-point BVP can be solved
with staggered FD method with exponential convergence at the end
points. The main idea was to link the  staggered FD approximation to
rational approximation of the exact NtD map and optimizing the
resulting approximation using Zolotorev theory. This method was
later extended to the solution of the two-sided problems by
splitting the solution into odd and even parts and solving two
one-sided problems on half-intervals using dual grids. Formerly
called optimal FD grids, the basic idea of spectrally matched FD
grids is summarized below.

Let us introduce the FD grid steps $\hat h_i,h_i$, $i=1,\ldots,k$.
We split the DtN map into odd and even parts and compute each of
them using a FD scheme  on  half interval. The odd and even problems
can respectively be written in mutually dual form as:
    \bea\label{1stfd}
            \frac {u^o_{i+1}-u^o_i}{h_i} = sv^o_i, \
            \frac {v^o_{i}-v^o_{i-1}}{\hat h_i}=su^o_i, \qquad i=1,\ldots,k, \quad u_{k+1}=0, \\
        \frac {u^e_{i+1}-u^e_i}{\hat h_i}=sv^e_i, \
        \frac {v^e_{i}-v^e_{i-1}}{h_i}=su^e_i, \qquad i=1,\ldots,k, \quad v_{k+1}=0. \nonumber
    \eea
It is known \cite{DrKn} that,
    \be \label{fraction} \frac{u_1^o} {v_1^o}=\frac{v_1^e}{u_1^e}=f_k(s)
        = \cfrac{1}{\hat h_1s+
            \cfrac{1}{h_1s+
            \cfrac{1}{\hat h_2s+\dots
            \cfrac{1}{h_{k-1}s+
            \cfrac{1}{\hat h_{k}s+
            \cfrac{1}{h_ks}}}}}}.
    \ee
Combining odd and even parts we obtain,
     \be\label{oddeven} u(0) = u_1^e+u_1^o, \
        u(1)=u_1^e-u_1^o,\quad
        v_1^e\approx \frac{v(0)+v(1)}{2}, \
        v_1^o\approx \frac{v(0)-v(1)}{2},
    \ee
and the FD-NtD as
    \be\label{impfd}\hat F
        = Z\left[\begin{array}{cc}f_k & 0 \\ 0 & 1/f_k\end{array}\right]Z^*.
    \ee

Construction of spectrally matched grids involves a reverse
procedure. First, rational approximation theory is used to obtain
$f_k$ that approximates the NtD map. The resulting rational function
is then used in (\ref{fraction}) to compute the grid steps $\hat
h_i,h_i$ using simple **** algorithm \cite{}. **** algorithm also
constructively shows that any [2k-1/2k] rational function can be
converted into an equivalent FD grid.

\subsection{Equivalence of discrete problems}
In this section, we show that if the number of finite elements are
chosen to be even ($n=2k$), the approximate NtD maps from FE and FD
grids are equivalent.
\begin{lemma}
For any set of parameters $l_i\in C$, $l=1,\ldots, 2k$ there exist
parameters $\hat h_i, h_i\in C\cup\infty$, $l=1,\ldots, k$, such
that
    \be\label{equiv}f(s)\equiv tanh(s/2)\ee
and vice versa.
\end{lemma}
\proof \ \ For any set of parameters $\hat h_i, h_i\in C\cup\infty$
there exist  polynomials $p$ and $q$ (at most) degree $k-1$ and $k$
respectively, such that the continued fraction expansion
(\ref{fraction}) can be presented as
    \[f_k=\frac{sp(s^2)}{q(s^2)},\]
and vice versa. Equating numerator and denominator of $f_k$ and
$tanh$ we equivalently transform (\ref{equiv}) to  polynomial
identities
    \be\label{polmatch} s p(s^2)\equiv t(s)-t(-s),\qquad q(s^2)\equiv t(s)+t(-s).\ee
Since $p$, $q$ and $t$ can be arbitrary polynomials of degree $k-1$,
$k$ and $2k$ respectively, then for any $p,q$ there is $t$
satisfying (\ref{polmatch}) and vice versa. {\bf (Aren't these
different p and q? If so, it is important not to confuse the
polynomials $p$ and $q$ with the polynomials in the second equation
in section 2. Should be rename these as $\tilde{p}$ and
$\tilde{q}$?)} {\color{blue} Agree, will do it later} \square

From the lemma, (\ref{impfe}) and (\ref{impfd}), we obtain the
following result about equivalence of the FE and FD DtN maps.
\begin{proposition}
If (\ref{equiv}) is valid, then
    \[\tilde F(s)\equiv \hat F(s).\]
\end{proposition}
Formula (\ref{polmatch}) can be used for computing the equivalent FE
from the FD and vice versa.

If the DtN maps are identical, then formula (\ref{propfe}) can also
be used for computing the propagator matrix for the FD
approximation.

If  $exp(s)$ matches $\exp(s)$ in $n$ non-coinciding frequencies,
then $f_k$  matches $\tanh(s)$ at the same frequencies, and $f_k$ is
Stieltjes function, $h_i,\hat h_i$ are real positive, and the
problem becomes Hermitian.

\section{Approximation of exterior problems}

Let a discretized  interval $\Omega_1=[x_-,x_+]$ have the propagator
matrix (from left to right),
\[{\cal \tilde G}=\left[\begin{array}{cc}
  exp_1(s) & 0 \\ 0 & exp_1(-s)
\end{array}\right]\]
in the spectral coordinates, where $exp_1(s)=t_1(-s)/t_1(s)$ defined
as in the previous section. First, let us impose the reflection
coefficient $h_2(s)/h_2(-s)$ at $x_+$, i.e., at the right boundary
any nontrivial solution can be represented as
$w(x_+)=[ch_2(s),-ch_2(-s)]^*$ in the spectral coordinates, where
$c\ne 0$ is an arbitrary constant. Then the reflection coefficient
at the left boundary will be the ratio of the components of
$w(x_-)={\cal \tilde G}^{-1} w(x_+)$.
 That is, \be\label{prodref}
  exp_1(-s)^2\frac{h_2(s)}{h_2(s)}=\frac{t_1(s)^2}{t_1(-s)^2}
\frac{h_2(s)}{h_2(-s)}.\ee If we impose the Dirichlet condition at
the right boundary of $\Omega_1$, which corresponds to $h_2(s)=1$,
then the reflection coefficient will be
$\frac{t_1(s)^2}{t_1(-s)^2}$. Let us now assume that we have a
connected interval $\Omega = \Omega_1 \cup \Omega_2$ with the
Dirichlet condition at the right boundary ($\Omega_2$ is assumed to
be on the right), and $\frac {h_2(s)} {h_2(-s)}$ is the reflection
coefficient of $\Omega_2$. Then (\ref{prodref}) would yield the
reflection coefficient of $\Omega$ that is just the product of the
reflection coefficients of the two subdomains.

Now, let as assume, that we use  the discrete problem in $\Omega$
for the approximation of (\ref{NtD1st}), i.e., $h(s)$ from
(\ref{Newman}) can be presented as
 $h(s)= t_1(s)^2h_2(s)$.
   If we set $h_1 \equiv t_1^2$ and $t_1 \equiv t_e$
and $t_2 \equiv t_p$, then the reflection coefficient of $\Omega$
will be identical to the one discussed in Section~2. However,
Dirichlet condition on the right of $\Omega_2$ makes it a one-sided
problem, and it is not necessary to restrict to the two-sided
approximation in the previous section. In fact, the original FD
optimal grids are optimized for the one-sided problems and can be
used effectively for $\Omega_2$. This approximation is equivalent to
the odd part of the FD approximation (\ref{1stfd}), i.e.,
\[\frac {u^o_{i+1}-u^o_i}{h_i} = sv^o_i, \
  \frac {v^o_{i}-v^o_{i-1}}{\hat h_i}=su^o_i, \qquad i=1,\ldots,k, \quad u_{k+1}=0.\]
Then $h_2$ can be obtained from the equality $f_k(s) = \frac
{h_2(s)-h_2(-s)} {h_2(s)+h_2(-s)}$, i.e., it can be an arbitrary
polynomial of degree $2k$.
$\cal  f $

\end{document}